\documentclass[11pt]{amsart}
\usepackage{amsfonts,amssymb,amsthm}
\usepackage{booktabs}
\usepackage[numbers,sort&compress]{natbib}
\usepackage[leqno]{amsmath}
\usepackage[mathscr]{euscript}
\usepackage{graphicx} 
\usepackage{subcaption}
\usepackage{booktabs}
\usepackage{color}
\usepackage{bm}
\usepackage{cases}
\usepackage{mathrsfs}
\usepackage{paralist}
\usepackage{hyperref}
\usepackage{framed}
\usepackage{epstopdf}
\usepackage[ruled, vlined]{algorithm2e}
\usepackage{caption}
\usepackage{tablefootnote}
\newtheorem {theorem}{Theorem}[section]
\newtheorem {proposition}{Proposition}[section]

\newtheorem {lemma}{Lemma}[section]
\newtheorem {example}{Example}[section]
\newtheorem {definition}{Definition}[section]

\makeatletter
\newcommand{\leqnomode}{\tagsleft@true}
\newcommand{\reqnomode}{\tagsleft@false}
\makeatother

\def\ees{{\accent"5E e}\kern-.385em\raise.2ex\hbox{\char'23}\kern-.08em}
\def\EES{{\accent"5E E}\kern-.5em\raise.8ex\hbox{\char'23 }}
\def\ow{o\kern-.42em\raise.82ex\hbox{
   \vrule width .12em height .0ex depth .075ex \kern-0.16em \char'56}\kern-.07em}
\def\OW{O\kern-.460em\raise1.36ex\hbox{
\vrule width .13em height .0ex depth .075ex \kern-0.16em \char'56}\kern-.07em}

\title{Optimality and Exact Moment-SDP Formulations for Robust Two-Stage Convex Polynomial Multiobjective Optimization}

\author{Yongqi Jiang}
\address[Yongqi Jiang]{Department of Mathematics, College of Sciences, Yanbian University, Yanji 133002, China}
\email{jiangyq0907@126.com}

\author{Zhe Hong}
\address[Zhe Hong]{Department of Mathematics, College of Sciences, Yanbian University, Yanji 133002, China}
\email{hongchul1990@gmail.com}

\author{Chengmiao Yang$^{\dag}$}
\address[Chengmiao Yang]{Academy for Advanced Interdisciplinary Studies, Northeast Normal University, Changchun 130024, Jilin Province, China}
\email{cmyang@nenu.edu.cn}

\thanks{The first two authors contribute equally.}

\thanks{$^{\dag}$Corresponding Author}



\date{\today}

\begin{document}

\maketitle

\begin{abstract} 
In this paper, we study Geoffrion properly efficient solutions of affinely adjustable robust two-stage convex polynomial multiobjective optimization problems with spectrahedral uncertainty.
For convex polynomial data, we establish sufficient optimality conditions and necessary conditions under a Slater condition.
Under SOS-convexity, we obtain an exact SOS characterization of proper efficiency.
Using strictly positive weighted-sum scalarization, we develop a finite-dimensional moment-SDP formulation and establish its optimal-value exactness under the strict feasibility condituon. 
We further establish the recovery results connecting Geoffrion properly efficient solutions with moment-SDP optimal solutions.

\end{abstract}


\section{Introduction}
Multiobjective optimization is an important topic of mathematical programming, with applications in resource allocation, production scheduling, and energy management, e.g., \cite{hao2024constrained,chankong2008multiobjective}.
However, model parameters in practical decision-making problems are often uncertain due to forecasting errors or changes in the external environment.
Therefore, ensuring feasibility and satisfactory performance under uncertainty has become a significant research issue. 
As a systematic approach to handling data uncertainty, robust optimization assumes that uncertain parameters take values in a prescribed uncertainty set and requires the constraints to hold for every realization in that set, e.g.,\cite{ben2009robust,bertsimas2011theory,chuong2020robust}. 
In the standard static robust optimization framework, all decision variables are determined before the uncertainty is observed, which can lead to conservative decisions.
Adjustable robust optimization addresses this limitation by allowing some decision variables to adapt to observed uncertainty, thereby providing a more flexible framework for sequential decision-making under uncertainty, e.g.,
\cite{ben2004adjustable,yanikouglu2019survey}.

Many practical problems involve decisions made at different stages. 
Some decisions must be made before the uncertain parameters are realized, while others can be adjusted after the relevant information becomes partially or fully available, e.g., \cite{ascione2016multi,li2015multi}. 
The two-stage optimization formulation captures this structure through here-and-now decisions and subsequent wait-and-see decisions \cite{ben2020tractable}. 
This motivates us to study two-stage multiobjective optimization problems under data uncertainty.

An uncertain convex polynomial multiobjective optimization problem is given by
\begin{align*}
\min\limits_{x\in\mathbb{R}^{n}} \quad
& \Bigl(f_1(x,u),\ldots,f_p(x,u)\Bigr)
\tag{$\rm UCP$}\label{UCP}\\
\text{s.t.} \quad \,
& g_j(x,u)\leq 0,\ j=1,\ldots,q,
\end{align*}
where $u\in U\subseteq\mathbb{R}^{s}$ denotes the uncertain
parameter, and $U$ is a nonempty compact uncertainty set. 
The objective functions $f_r:\mathbb{R}^{n}\times\mathbb{R}^{s}\rightarrow
\mathbb{R}$, $r=1,\ldots,p$, and the constraint functions
$g_j:\mathbb{R}^{n}\times\mathbb{R}^{s}\rightarrow\mathbb{R}$,
$j=1,\ldots,q$, are given by
\begin{equation*}
f_r(x,u):=f_r^0(x)+\sum_{\ell=1}^{s}u_\ell f_r^\ell(x),
\quad
g_j(x,u):=g_j^0(x)+\sum_{\ell=1}^{s}u_\ell g_j^\ell(x),
\end{equation*}
where $f_r^\ell,g_j^\ell,$
$\ell=0,1,\ldots,s$ are given real polynomials. 
For each fixed $u\in U$, the functions $f_r(\cdot,u)$ and $g_j(\cdot,u)$ are assumed to be convex on $\mathbb{R}^n$.

Extending problem \eqref{UCP} to a two-stage setting, we consider the following uncertain two-stage convex polynomial multiobjective optimization problem:
\begin{align*}
\min\limits_{x\in\mathbb{R}^{n},\,y(\cdot)}\quad \,
& \Bigl(f_1(x,u)+v_1^{\top}y(u),\ldots,
f_p(x,u)+v_p^{\top}y(u)\Bigr)
\tag{$\rm UTCP$}\label{UTCP}\\
\text{s.t.}\qquad 
& g_j(x,u)+\theta_j^{\top}y(u)\leq 0,
\ j=1,\ldots,q,\ \forall u\in U.
\end{align*}
Here, $x\in\mathbb{R}^{n}$ is the first-stage decision variable, and $y(\cdot):U\rightarrow\mathbb{R}^{m}$ is the second-stage decision rule depending on the uncertain parameter $u\in U$. Let $v_r\in\mathbb{R}^m,\, r=1,\ldots,p$ and $\theta_j\in\mathbb{R}^m,\, j=1,\ldots,q$ be given vectors representing the coefficients of the second-stage decision variable in the $r$-th objective function and the $j$-th constraint function, respectively.

Since the second-stage decision rule $y(\cdot)$ is an arbitrary function of the uncertain parameter, problem \eqref{UTCP} is generally infinite-dimensional and difficult to solve directly. To facilitate a tractable  formulation, we assume that the uncertainty set $U$ is a nonempty bounded spectrahedron of the form (see e.g., \cite{ramana1995some,vinzant2014spectrahedron})
\[
U:=\left\{
u=(u_1,\ldots,u_s)\in\mathbb{R}^{s}
\ \middle|\
H_0+\sum_{\ell=1}^{s}u_\ell H_\ell\succeq 0
\right\},
\]
where $H_\ell$, $\ell=0,1,\ldots,s$, are given
symmetric matrices. 
We further restrict the second-stage rule to be affine (see e.g., \cite{chen2009uncertain,chuong2020generalized}):
\[
y(u):=y^0+Wu,\ u\in U,
\]
where $y^0\in\mathbb{R}^{m}$ and $W\in\mathbb{R}^{m\times s}$.

By imposing robust feasibility over $U$ and taking the worst-case value of each objective function, we obtain the adjustable robust counterpart of problem
\eqref{UTCP} as follows:
\begin{align*}
\min\limits_{\substack{
x\in\mathbb{R}^{n},\,y^0\in\mathbb{R}^{m},\\
W\in\mathbb{R}^{m\times s}}}
\quad \,\,
& \Biggl(\max_{u\in U}\Bigl\{f_1(x,u)+v_1^\top y(u)\Bigr\},\ldots,\max_{u\in U}\Bigl\{f_p(x,u)+v_p^\top y(u)\Bigr\}\Biggr)\tag{RTCP}\label{RTCP}\\
\text{s.t.} \quad \qquad 
& g_j(x,u)+\theta_j^\top (y^0+Wu)\leq 0,
\ j=1,\ldots,q,\ \forall u\in U.
\end{align*}
Problem~\eqref{RTCP} is the main model studied in this paper.

In multiobjective optimization, efficient solutions are commonly used to describe the trade-offs among conflicting objectives. However, efficiency alone does not exclude unbounded trade-offs, where the ratio of improvement in one objective to deterioration in another can become arbitrarily large. Proper efficiency was introduced to exclude such unbounded trade-offs. 
In this paper, we focus on Geoffrion properly efficient solutions of problem~\eqref{RTCP}, and use Benson proper efficiency as an equivalent geometric characterization, e.g., \cite{benson1978existence,benson1977vector,geoffrion1968proper,white1983concepts}. The formal definitions of these solution concepts are given in Sect.~\ref{sect:2}.

Existing studies have developed optimality conditions and numerical methods for adjustable robust multiobjective optimization with linear, quadratic, and general convex data, e.g., \cite{chuong2025adjustable,chuong2024solving,nguyen2026efficiency}. 
Related studies have investigated robust optimality and solution existence using minimax and tangency-based methods \cite{hong2022minimax,jiao2026robust}. 
In parallel, optimality characterizations and semidefinite programming approaches have been developed for robust multiobjective optimization with SOS-convex polynomial data, particularly in single-stage models, e.g., \cite{huang2025robust,liguo2021finding}. 
However, extending these approaches to affinely adjustable robust two-stage convex polynomial multiobjective problems remains challenging. 
This motivates us to combine convex optimality theory with SOS techniques and semidefinite programming to characterize and compute Geoffrion properly efficient solutions for the proposed model.

\subsection*{Contributions}
In this paper, we study a class of affinely adjustable robust two-stage convex polynomial multiobjective optimization problems under spectrahedral uncertainty, focusing on the optimality characterization of Geoffrion properly efficient solutions and their numerical computation. 
The main contributions of this paper are as follows.
\begin{itemize}
\item[{\rm (C1)}] 
For convex polynomial data, we establish necessary and sufficient optimality conditions for Geoffrion properly efficient solutions of problem \eqref{RTCP}. 
The necessity result requires the Slater condition, whereas the sufficiency result does not. 
\item[{\rm (C2)}] 
Under SOS-convexity, we establish an exact SOS characterization of
Geoffrion proper efficiency by replacing the global polynomial
nonnegativity condition in the optimality system with an SOS
constraint.
\item[{\rm (C3)}] 
Based on weighted-sum scalarization with strictly positive weights, we construct a finite-dimensional moment-SDP relaxation for the scalarized problem.
We establish its optimal-value exactness and solution-recovery properties. A numerical example is also provided to illustrate the effectiveness of the proposed approach.
\end{itemize}

\subsection*{Outline}
The remainder of this paper is organized as follows.
In Sect.~\ref{sect:2}, we recall some basic notion and and establish a preliminary result on weighted-sum scalarization. 
In Sect.~\ref{sect:3}, we derive necessary and sufficient optimality conditions for Geoffrion properly efficient solutions with convex polynomial data and examine the role of the Slater condition. In Sect.~\ref{sect:4}, we specialize the analysis to SOS-convex polynomial data, establish an exact SOS characterization of proper efficiency, develop the corresponding moment-SDP formulation and solution-recovery results, and present a numerical example. Conclusions are given in Sect.~\ref{sect:5}.

\section{Preliminaries}\label{sect:2}
In this section, we recall some notation and preliminary concepts used in the subsequent analysis. Let $\mathbb{N}_{+}:=\{1,2,\ldots\}$. For any $d\in\mathbb{N}_{+}$, $\mathbb{R}^{d}$ denotes the $d$-dimensional Euclidean space equipped with the inner product $\langle x,z\rangle:=x^{\top}z$ and the induced norm $\|x\|:=\sqrt{\langle x,x\rangle}$. We denote by $\mathbb{S}^{k}$ the space of real symmetric $k\times k$ matrices and by $\mathbb{S}^{k}_{+}$ the cone of positive semidefinite matrices. For $A\in\mathbb{S}^{k}$, $A\succeq 0$ and $A\succ 0$ mean that $A$ is positive semidefinite and positive definite, respectively. For matrices $X$ and $Y$ of the same size, their inner product is defined by $\langle X,Y\rangle:=\operatorname{tr}(X^{\top}Y)$. Let $\mathbb{R}^{p}_{+}$ and $\mathbb{R}^{p}_{++}$ denote the nonnegative and strictly positive orthants of $\mathbb{R}^{p}$. For $z,w\in\mathbb{R}^{p}$, we write $z\leq w$ if $w-z\in\mathbb{R}^{p}_{+}$ and $z<w$ if $w-z\in\mathbb{R}^{p}_{++}$. For a set $C$, $\operatorname{cl}C$ denotes its closure and $\operatorname{cone}C$ denotes its conic hull.

For later use, define the robust objective and constraint functions associated with problem~\eqref{RTCP} by
\[
F_r(x,y^0,W)
:=
\max_{u\in U}
\left\{
f_r(x,u)+v_r^{\top}(y^0+Wu)
\right\},
\ r=1,\ldots,p,
\]
\[
G_j(x,y^0,W)
:=
\max_{u\in U}
\left\{
g_j(x,u)+\theta_j^{\top}(y^0+Wu)
\right\},
\ j=1,\ldots,q.
\]
Since $U$ is compact and the functions are continuous in $u$, the maxima defining $F_r$ and $G_j$ are attained. Moreover, each $F_r$ and $G_j$ is convex in $(x,y^0,W)$. Accordingly, the robust vector-valued objective mapping is defined as
\[
F(x,y^0,W)
:=
\bigl(
F_1(x,y^0,W),\ldots,F_p(x,y^0,W)
\bigr),
\]
so the feasible set of problem \eqref{RTCP} can be equivalently represented as
\[
\mathcal F
:=
\left\{
(x,y^0,W)
\in
\mathbb R^n\times\mathbb R^m\times\mathbb R^{m\times s}
\;\middle|\;
G_j(x,y^0,W)\leq0,\ j=1,\ldots,q
\right\}.
\]

Based on the robust objective mapping $F(x,y^0,W)$ and the feasible set $\mathcal F$, we now introduce the notions of efficient and Geoffrion properly efficient solutions for problem \eqref{RTCP}.

\begin{definition}[Efficient solution]\label{def:efficient}
A feasible solution $(x^*,y^{0*},W^*)\in\mathcal{F}$ is called an efficient solution of problem \eqref{RTCP} if there does not exist any $(x,y^0,W)\in\mathcal{F}$ such that
\[
F_r(x,y^0,W)\leq F_r(x^*,y^{0*},W^*),
\ r=1,\ldots,p,
\]
\[
F_{r_0}(x,y^0,W)<F_{r_0}(x^*,y^{0*},W^*),
\
\text{for some } r_0\in\{1,\ldots,p\}.
\]
\end{definition}

\begin{definition}[Geoffrion properly efficient solution]
\label{def:Geoffrion-proper}
A feasible solution $(x^*,y^{0*},W^*)\in\mathcal{F}$ is called a
Geoffrion properly efficient solution of problem \eqref{RTCP} if it is efficient and there exists a real number $M>0$ such that, for all $r\in\{1,\ldots,p\}$ and $(x,y^0,W)\in\mathcal{F}$ satisfying $F_r(x,y^0,W)
<
F_r(x^*,y^{0*},W^*),$ there exists an index
$k\in\{1,\ldots,p\}$ such that $F_k(x^*,y^{0*},W^*)
<
F_k(x,y^0,W)$ and
\[
\frac{
F_r(x^*,y^{0*},W^*)-F_r(x,y^0,W)
}{
F_k(x,y^0,W)-F_k(x^*,y^{0*},W^*)
}
\leq M.
\]
\end{definition}

To provide a geometric characterization of properly efficient solutions in the objective space, define the objective image set of problem \eqref{RTCP} by
\[
\mathcal Y
:=
F(\mathcal F)
=
\left\{
F(x,y^0,W)
\;\middle|\;
(x,y^0,W)\in\mathcal F
\right\}
\subseteq \mathbb R^p.
\]
Since $\mathcal F$ is convex and $F$ is componentwise convex, the upper image $\mathcal Y+\mathbb R_+^p$ is convex. We next introduce the notion of a Benson properly efficient solution for problem \eqref{RTCP}.

\begin{definition}[Benson properly efficient solution]
\label{def:Benson-proper}
A feasible solution $(x^*,y^{0*},W^*)\in\mathcal F$ is called a Benson properly efficient solution of problem \eqref{RTCP} if
\[
\operatorname{cl}\operatorname{cone}
\left(
\mathcal Y+\mathbb R_+^p-F(x^*,y^{0*},W^*)
\right)
\cap
\left(-\mathbb R_+^p\right)
=
\{0\}.
\]
\end{definition}

By \cite[Theorem~2.48]{ehrgott2005multicriteria}, Geoffrion proper efficiency and Benson proper efficiency are equivalent for problem \eqref{RTCP}. 
Since the upper image  $\mathcal{Y}+\mathbb{R}_{+}^{p}$ is convex, the geometric characterization of Benson proper efficiency leads to the following weighted-sum scalarization result.
%
%
\begin{lemma}
\label{lem:positive-supporting-vector}
Let $(x^*,y^{0*},W^*)\in\mathcal F$ be a Benson properly efficient solution of problem \eqref{RTCP}. Then there exists a vector $\alpha=(\alpha_1,\ldots,\alpha_p)\in\mathbb R_{++}^p$, such that
\[
\sum_{r=1}^{p}\alpha_r F_r(x,y^0,W)
\geq
\sum_{r=1}^{p}\alpha_r F_r(x^*,y^{0*},W^*),
\
\forall (x,y^0,W)\in\mathcal F.
\]
\end{lemma}

\begin{proof}
Define
\[
\mathcal C
:=
\operatorname{cl}\operatorname{cone}
\left(
\mathcal Y+\mathbb R_+^p-F(x^*,y^{0*},W^*)
\right).
\]
Since the upper image $\mathcal Y+\mathbb R_+^p$ is convex, $\mathcal C$ is a closed convex cone. Moreover, the Benson proper efficiency of $(x^*,y^{0*},W^*)$ implies that $\mathcal C\cap\left(-\mathbb R_+^p\right)=\{0\}$.
Let
\[
\mathcal C^*
:=
\left\{
\mu\in\mathbb R^p
\;\middle|\;
\mu^\top d\geq0
\text{ for all }d\in\mathcal C
\right\}
\]
be the dual cone of $\mathcal C$. We first show that
\[
\mathcal C^*\cap\mathbb R_{++}^p\neq\varnothing.
\]

Suppose, to the contrary, that $\mathcal C^*\cap\mathbb R_{++}^p=\varnothing$. Since $\mathcal C^*$ is a nonempty closed convex cone and $\mathbb R_{++}^p$ is a nonempty open convex set, the separation theorem yields a nonzero vector $\gamma\in\mathbb R^p$ and a scalar $\beta\in\mathbb R$ such that
\begin{equation}\label{eq:proper-separation}
\gamma^\top\mu
\leq
\beta
\leq
\gamma^\top\nu,
\
\forall\mu\in\mathbb R_{++}^p,\
\forall\nu\in\mathcal C^*.
\end{equation}
Since $0\in\mathcal C^*$,  \eqref{eq:proper-separation} implies $\beta\leq0$. 
For any fixed $\mu\in\mathbb R_{++}^p$ and $t>0$, we have $t\mu \in R_{++}^p$. 
Letting $t\to0^+$ gives $\beta\geq0$. 
Thus, $\beta=0$, and consequently,
\[
\gamma^\top\mu\leq0,
\quad
\forall\mu\in\mathbb R_{++}^p.
\]
For each $r=1,\ldots,p$, let $e_r$ denote the $r$-th standard basis vector of $\mathbb R^p$. Taking
$\mu^\varepsilon:=e_r+\varepsilon e$, where
$e:=(1,\ldots,1)^\top$ and $\varepsilon>0$, we have
$\mu^\varepsilon\in\mathbb R_{++}^p$ and hence
$\gamma^\top\mu^\varepsilon\leq0$. Letting $\varepsilon\to0^+$ yields
$\gamma_r=\gamma^\top e_r\leq0$. Since $r$ is arbitrary, $\gamma\in-\mathbb R_+^p$.
Moreover, since $\beta=0$, the right-hand inequality in
\eqref{eq:proper-separation} yields $\gamma^\top\nu\geq0,
\
\forall\nu\in\mathcal C^*$, and hence $\gamma\in\mathcal C^{**}$. Since $\mathcal C$ is a closed convex cone,
the bipolar theorem gives $\mathcal C^{**}=\mathcal C$, and therefore
$\gamma\in\mathcal C$. Consequently,
\[
0\neq\gamma\in
\mathcal C\cap\left(-\mathbb R_+^p\right),
\]
which contradicts the Benson proper efficiency condition. It follows that
\[
\mathcal C^*\cap\mathbb R_{++}^p\neq\varnothing.
\]

We choose $\alpha=(\alpha_1,\ldots,\alpha_p)
\in
\mathcal C^*\cap\mathbb R_{++}^p$. For any $(x,y^0,W)\in\mathcal F$, the definition $\mathcal Y=F(\mathcal F)$ implies that $F(x,y^0,W)\in\mathcal Y$. Since $0\in\mathbb R_+^p$, it follows that
\[
F(x,y^0,W)-F(x^*,y^{0*},W^*)
\in
\mathcal Y+\mathbb R_+^p-F(x^*,y^{0*},W^*)
\subseteq\mathcal C.
\]
Since $\alpha\in\mathcal C^*$, we obtain
\begin{equation*}
0\leq\alpha^\top
\left(F(x,y^0,W)-F(x^*,y^{0*},W^*)\right)=\sum_{r=1}^{p}\alpha_rF_r(x,y^0,W)-\sum_{r=1}^{p}\alpha_rF_r(x^*,y^{0*},W^*).
\end{equation*}
Therefore,
\[
\sum_{r=1}^{p}\alpha_rF_r(x,y^0,W)
\geq
\sum_{r=1}^{p}\alpha_rF_r(x^*,y^{0*},W^*),
\quad
\forall(x,y^0,W)\in\mathcal F.
\qedhere
\]
\end{proof}

\section{Optimality for Two-stage Convex Polynomial Multiobjective Problems}\label{sect:3}
This section establishes necessary and sufficient optimality conditions for Geoffrion properly efficient solutions of problem~\eqref{RTCP} with convex polynomial data. 

We say that the Slater condition holds if there exists $(\bar{x},\bar{y}^0,\bar{W})
\in
\mathbb R^n\times\mathbb R^m\times\mathbb R^{m\times s}$, such that
\begin{equation}\label{eq:slater-condition}
g_j(\bar{x},u)
+
\theta_j^\top\left(\bar{y}^0+\bar{W}u\right)
<0,
\
\forall u\in U,\ j=1,\ldots,q.
\end{equation}
For notational convenience, we identify
$\mathbb{R}^{n}\times\mathbb{R}^{m}\times\mathbb{R}^{m\times s}$
with $\mathbb{R}^{n+m+ms}$ throughout this section.

\begin{theorem}[Necessary optimality conditions]
\label{thm:necessary-optimality}
Let $(x^*,y^{0*},W^*)\in\mathcal{F}$ be a Geoffrion properly efficient solution of problem \eqref{RTCP}. Suppose that the Slater condition \eqref{eq:slater-condition} holds. Then there exist $\alpha=(\alpha_1,\ldots,\alpha_p)\in\mathbb{R}_{++}^p,
\
\alpha_r^\ell\in\mathbb{R},
\ r=1,\ldots,p,\ \ell=1,\ldots,s$,
and
$\lambda=(\lambda_1,\ldots,\lambda_q)\in\mathbb{R}_+^q,
\
\lambda_j^\ell\in\mathbb{R},
\ j=1,\ldots,q,\ \ell=1,\ldots,s,$
such that
\begin{equation}\label{eq:necessary-lmi-conditions}
\alpha_rH_0+\sum_{\ell=1}^{s}\alpha_r^\ell H_\ell\succeq0,
\ r=1,\ldots,p,
\quad
\lambda_jH_0+\sum_{\ell=1}^{s}\lambda_j^\ell H_\ell\succeq0,
\ j=1,\ldots,q,
\end{equation}
\begin{equation}\label{eq:necessary-balance-conditions}
\sum_{r=1}^{p}\alpha_rv_r+\sum_{j=1}^{q}\lambda_j\theta_j=0,
\quad
\sum_{r=1}^{p}\alpha_r^\ell v_r
+\sum_{j=1}^{q}\lambda_j^\ell\theta_j=0,
\ \ell=1,\ldots,s,
\end{equation}

\begin{equation}\label{eq:necessary-polynomial-condition}
\begin{aligned}
\Phi(x)
:={}&
\sum_{r=1}^{p}
\left(
\alpha_rf_r^0(x)
+\sum_{\ell=1}^{s}\alpha_r^\ell f_r^\ell(x)
\right)
+\sum_{j=1}^{q}
\left(
\lambda_jg_j^0(x)
+\sum_{\ell=1}^{s}\lambda_j^\ell g_j^\ell(x)
\right) \\
-&
\sum_{r=1}^{p}\alpha_rF_r(x^*,y^{0*},W^*).\\
\Phi(x) \geq &\
\Phi(x^*)=0,
\
\forall x\in\mathbb{R}^n.
\end{aligned}
\end{equation}
\end{theorem}

\begin{proof}
By the equivalence of Geoffrion and Benson proper efficiency and Lemma~\ref{lem:positive-supporting-vector}, there exists   $\alpha=(\alpha_1,\ldots,\alpha_p)\in\mathbb R_{++}^p$, such that
\[
\sum_{r=1}^{p}\alpha_rF_r(x,y^0,W)
\geq
\sum_{r=1}^{p}\alpha_rF_r(x^*,y^{0*},W^*),
\quad
\forall (x,y^0,W)\in\mathcal F.
\]
It follows that the following system of convex inequalities has no solution:
\begin{equation}\label{eq:infeasible-convex-system}
\left\{
\begin{aligned}
&\sum_{r=1}^{p}\alpha_r
\left[
F_r(x,y^0,W)-F_r(x^*,y^{0*},W^*)
\right]<0,\\
&G_j(x,y^0,W)<0,
\quad j=1,\ldots,q.
\end{aligned}
\right.
\end{equation}
Since all the functions in \eqref{eq:infeasible-convex-system} are finite-valued and convex on $\mathbb R^{n+m+ms}$, Rockafellar's convex alternative theorem \cite[Theorem~21.1]{rockafellar2015convex} yields $\tau\geq0$ and $\eta_j\geq0$, $j=1,\ldots,q$, not all zero, such that
\begin{equation}\label{eq:alternative-inequality}
\tau\sum_{r=1}^{p}\alpha_r
\left[
F_r(x,y^0,W)-F_r(x^*,y^{0*},W^*)
\right]
+\sum_{j=1}^{q}\eta_jG_j(x,y^0,W)
\geq0,
\
\forall (x,y^0,W)\in\mathbb R^{n+m+ms}.
\end{equation}

We next show that $\tau>0$. If $\tau=0$, then there exists
$j_0\in\{1,\ldots,q\}$ such that $\eta_{j_0}>0$. By the Slater condition \eqref{eq:slater-condition}, we have
\[
\sum_{j=1}^{q}
\eta_jG_j(\bar{x},\bar{y}^0,\bar{W})<0,
\]
which contradicts \eqref{eq:alternative-inequality}. Therefore, $\tau>0$. Dividing \eqref{eq:alternative-inequality} by $\tau$ and setting $\lambda_j:=\eta_j/\tau,$ $j=1,\ldots,q$, we obtain
\begin{equation}\label{eq:normalized-alternative}
\sum_{r=1}^{p}\alpha_r
\left[
F_r(x,y^0,W)-F_r(x^*,y^{0*},W^*)
\right]
+
\sum_{j=1}^{q}\lambda_jG_j(x,y^0,W)
\geq0,
\
\forall (x,y^0,W)\in\mathbb R^{n+m+ms}.
\end{equation}
By the definitions of $F_r(x,y^0,W)$ and $G_j(x,y^0,W)$, \eqref{eq:normalized-alternative} can be rewritten as
\begin{equation}\label{eq:expanded-inf-max}
\begin{aligned}
\inf_{(x,y^0,W)\in\mathbb R^{n+m+ms}}
\Bigg\{&
\sum_{r=1}^{p}\alpha_r
\max_{u_1^r\in U}
\left[
f_r(x,u_1^r)
+
v_r^\top\left(y^0+Wu_1^r\right)
\right] \\
 + &
\sum_{j=1}^{q}\lambda_j
\max_{u_2^j\in U}
\left[
g_j(x,u_2^j)
+
\theta_j^\top\left(y^0+Wu_2^j\right)
\right]
\Bigg\}
\geq
\sum_{r=1}^{p}\alpha_rF_r(x^*,y^{0*},W^*).
\end{aligned}
\end{equation}
Let $\Omega:=U^p\times U^q$ and $\boldsymbol{u}
:=
\left(
u_1^1,\ldots,u_1^p,
u_2^1,\ldots,u_2^q
\right)
\in\Omega$.
Define $\Psi:\mathbb R^{n+m+ms}\times\Omega\rightarrow\mathbb R$ by
\begin{equation*}
\Psi(x,y^0,W,\boldsymbol{u})
:=\sum_{r=1}^{p}\alpha_r
\left[
f_r(x,u_1^r)+
v_r^\top\left(y^0+Wu_1^r\right)
\right] +
\sum_{j=1}^{q}\lambda_j
\left[
g_j(x,u_2^j)
+
\theta_j^\top\left(y^0+Wu_2^j\right)
\right].
\end{equation*}
Therefore,
\eqref{eq:expanded-inf-max} is equivalent to
\[
\inf_{(x,y^0,W)\in\mathbb R^{n+m+ms}}\
\max_{\boldsymbol{u}\in\Omega}
\Psi(x,y^0,W,\boldsymbol{u})
\geq
\sum_{r=1}^{p}\alpha_rF_r(x^*,y^{0*},W^*).
\]
Since $\Omega$ is convex, $\Psi$ is continuous, convex in $(x,y^0,W)$, and affine in $u$, the minimax theorem \cite[Theorem~4.2]{sion1958general} gives
\begin{equation}\label{eq:minimax-psi}
\begin{aligned}
\max_{\boldsymbol{u}\in\Omega}\
\inf_{(x,y^0,W)\in\mathbb R^{n+m+ms}}
\Psi(x,y^0,W,\boldsymbol{u})\
& =\
\inf_{(x,y^0,W)\in\mathbb R^{n+m+ms}}\
\max_{\boldsymbol{u}\in\Omega}
\Psi(x,y^0,W,\boldsymbol{u}) \\
& \geq \
\sum_{r=1}^{p}\alpha_rF_r(x^*,y^{0*},W^*).
\end{aligned}
\end{equation}
Since the mapping
\[
\boldsymbol{u}
\longmapsto
\inf_{(x,y^0,W)\in\mathbb R^{n+m+ms}}
\Psi(x,y^0,W,\boldsymbol{u})
\]
is upper semicontinuous on $\Omega$, the maximum in \eqref{eq:minimax-psi} is attained. Hence, there exists $\boldsymbol{u}^*
:=
\left(
u_1^{1*},\ldots,u_1^{p*},
u_2^{1*},\ldots,u_2^{q*}
\right)
\in\Omega$, such that
\[
\inf_{(x,y^0,W)\in\mathbb R^{n+m+ms}}
\Psi(x,y^0,W,\boldsymbol{u}^*)
\geq
\sum_{r=1}^{p}\alpha_rF_r(x^*,y^{0*},W^*).
\]
Equivalently,
\begin{equation}\label{eq:fixed-uncertainty}
\Psi(x,y^0,W,\boldsymbol{u}^*)
\geq
\sum_{r=1}^{p}\alpha_rF_r(x^*,y^{0*},W^*),
\
\forall(x,y^0,W)\in\mathbb R^{n+m+ms}.
\end{equation}

For each $r=1,\ldots,p$ and $j=1,\ldots,q$, let $u_1^{r*}
=
\left(
u_{1,1}^{r*},\ldots,u_{1,s}^{r*}
\right)\in U,
\
u_2^{j*}
=
\left(
u_{2,1}^{j*},\ldots,u_{2,s}^{j*}
\right)\in U$
and define $\alpha_r^\ell
:=
\alpha_r u_{1,\ell}^{r*},
$
and $\lambda_j^\ell
:=
\lambda_j u_{2,\ell}^{j*},
\
\ell=1,\ldots,s.$
Since $u_1^{r*}\in U$, it follows that
\[
H_0+\sum_{\ell=1}^{s}u_{1,\ell}^{r*}H_\ell\succeq0.
\]
Multiplying this matrix inequality by $\alpha_r>0$, we obtain
\[
\alpha_rH_0+\sum_{\ell=1}^{s}\alpha_r^\ell H_\ell
=
\alpha_r
\left(
H_0+\sum_{\ell=1}^{s}u_{1,\ell}^{r*}H_\ell
\right)
\succeq0,
\ r=1,\ldots,p.
\]
Similarly, since $u_2^{j*}\in U$ and $\lambda_j\geq0$, we obtain
\[
\lambda_jH_0+\sum_{\ell=1}^{s}\lambda_j^\ell H_\ell
=
\lambda_j
\left(
H_0+\sum_{\ell=1}^{s}u_{2,\ell}^{j*}H_\ell
\right)
\succeq0,
\ j=1,\ldots,q.
\]
Let $W_\ell\in\mathbb R^m$ denote the $\ell$-th column of $W$,
$\ell=1,\ldots,s$. Using $Wu=\sum_{\ell=1}^{s}u_\ell W_\ell$ together with the affine dependence of $f_r$ and $g_j$ on $u$, inequality \eqref{eq:fixed-uncertainty} can be written as
\begin{equation}\label{eq:expanded-fixed-uncertainty}
\begin{aligned}
&\sum_{r=1}^{p}
\left(
\alpha_r f_r^0(x)
+\sum_{\ell=1}^{s}\alpha_r^\ell f_r^\ell(x)
\right)
+\sum_{j=1}^{q}
\left(
\lambda_j g_j^0(x)
+\sum_{\ell=1}^{s}\lambda_j^\ell g_j^\ell(x)
\right)\\
&\quad
-\sum_{r=1}^{p}\alpha_rF_r(x^*,y^{0*},W^*)
+\left(
\sum_{r=1}^{p}\alpha_rv_r
+\sum_{j=1}^{q}\lambda_j\theta_j
\right)^\top y^0\\
&\quad
+\sum_{\ell=1}^{s}
\left(
\sum_{r=1}^{p}\alpha_r^\ell v_r
+\sum_{j=1}^{q}\lambda_j^\ell\theta_j
\right)^\top W_\ell
\geq0,
\
\forall(x,y^0,W)\in\mathbb R^{n+m+ms}.
\end{aligned}
\end{equation}
Since \eqref{eq:expanded-fixed-uncertainty} holds for all
$(x,y^0,W)\in\mathbb R^{n+m+ms}$, and both $y^0$ and the columns
$W_\ell$, $\ell=1,\ldots,s$, are unrestricted, the coefficients of their corresponding linear terms must vanish. Hence,
\begin{equation}\label{eq:balance-conditions}
\sum_{r=1}^{p}\alpha_rv_r
+\sum_{j=1}^{q}\lambda_j\theta_j
=0,
\quad
\sum_{r=1}^{p}\alpha_r^\ell v_r
+\sum_{j=1}^{q}\lambda_j^\ell\theta_j
=0,
\
\ell=1,\ldots,s.
\end{equation}

Define the polynomial $\Phi:\mathbb R^n\rightarrow\mathbb R$ by
\begin{equation*}
\Phi(x):=\sum_{r=1}^{p}
\left(
\alpha_r f_r^0(x)
+\sum_{\ell=1}^{s}\alpha_r^\ell f_r^\ell(x)
\right) +\sum_{j=1}^{q}
\left(\lambda_j g_j^0(x)+\sum_{\ell=1}^{s}\lambda_j^\ell g_j^\ell(x)
\right)-\sum_{r=1}^{p}
\alpha_rF_r(x^*,y^{0*},W^*).
\end{equation*}
By \eqref{eq:expanded-fixed-uncertainty} and the balance conditions
\eqref{eq:balance-conditions}, we obtain
\begin{equation}\label{eq:global-nonnegativity}
\Phi(x)\geq0,
\
\forall x\in\mathbb R^n.
\end{equation}
We finally show that $\Phi(x^*)=0$. 
By \eqref{eq:fixed-uncertainty}, we have
\begin{equation*}
\Psi\left(x^{*}, y^{0*}, W^{*}, \boldsymbol{u}^{*}\right) \geq \sum_{r=1}^{p} \alpha_{r} F_{r}\left(x^{*}, y^{0 *}, W^{*}\right) .
\end{equation*}
On the other hand, by the definition of $F_{r}$, and the nonnegativity of $\alpha_{r}$  and $\lambda_{j}$, we obtain
\begin{equation*}
\Psi\left(x^{*}, y^{0*}, W^{*}, \boldsymbol{u}^{*}\right) \leq \sum_{r=1}^{p} \alpha_{r} F_{r}\left(x^{*}, y^{0 *}, W^{*}\right) .
\end{equation*}
Hence,
\begin{equation*}
\Psi\left(x^{*}, y^{0 *}, W^{*}, \boldsymbol{u}^{*}\right)=\sum_{r=1}^{p} \alpha_{r} F_{r}\left(x^{*}, y^{0 *}, W^{*}\right) .
\end{equation*}
By \eqref{eq:expanded-fixed-uncertainty},
\eqref{eq:balance-conditions}, and the definition of $\Phi$, this implies that  $\Phi\left(x^{*}\right)=0$. Combining this with \eqref{eq:global-nonnegativity}, we obtain $\Phi(x) \geq \Phi\left(x^{*}\right)=0, \, \forall x \in \mathbb{R}^{n}$.
This completes the proof.
%
%
%
%
%
\end{proof}

The multiplier conditions derived above are also sufficient for Geoffrion proper efficiency, without requiring the Slater condition.

\begin{theorem}[Sufficient optimality conditions]
\label{thm:sufficient-optimality}
Let $(x^*,y^{0*},W^*)\in\mathcal{F}$. 
Suppose that there exist  $\alpha \in \mathbb{R}_{++}^{p}, \ \lambda \in \mathbb{R}_{+}^{q} $, and $\alpha_{r}^{\ell},\ \lambda_{j}^{\ell} \in \mathbb{R},\ r=1, \ldots, p,\ j=1, \ldots, q,\ \ell=1, \ldots, s$, satisfying conditions \eqref{eq:necessary-lmi-conditions}--\eqref{eq:necessary-polynomial-condition} of Theorem \ref{thm:necessary-optimality}.
Then
$(x^*,y^{0*},W^*)$ is a Geoffrion properly efficient solution of problem
\eqref{RTCP}.
\end{theorem}

\begin{proof}
For each $r=1,\ldots,p$, define
\[u_1^r
:=
\left(
\frac{\alpha_r^1}{\alpha_r},
\ldots,
\frac{\alpha_r^s}{\alpha_r}
\right).
\]
Since $\alpha_r>0$, dividing the first matrix inequality in
\eqref{eq:necessary-lmi-conditions} by $\alpha_r$ gives
\[
H_0+\sum_{\ell=1}^{s}
\frac{\alpha_r^\ell}{\alpha_r}H_\ell
\succeq0.
\]
Hence, by the definition of $U$, $u_1^r\in U$, and
$\alpha_r^\ell=\alpha_ru_{1,\ell}^r$ for $\ell=1,\ldots,s$.
We next show that, for each $j=1,\ldots,q$, there exists
$u_2^j\in U$ satisfying $\lambda_j^\ell=\lambda_ju_{2,\ell}^j,
\
\ell=1,\ldots,s$. 
Similarly, if $\lambda_j>0$, define
\[
u_2^j
:=
\left(
\frac{\lambda_j^1}{\lambda_j},
\ldots,
\frac{\lambda_j^s}{\lambda_j}
\right).
\]
Dividing the second matrix inequality in \eqref{eq:necessary-lmi-conditions} by $\lambda_j$ yields $u_2^j\in U,
\
\lambda_j^\ell=\lambda_ju_{2,\ell}^j,
\
\ell=1,\ldots,s.$
Suppose now that $\lambda_j=0$. The second matrix inequality in \eqref{eq:necessary-lmi-conditions} becomes
\[
\sum_{\ell=1}^{s}\lambda_j^\ell H_\ell\succeq0.
\]
We claim that $\lambda_j^\ell=0$ for all $\ell=1,\ldots,s$. 
Otherwise, let $d_j:=
\left(
\lambda_j^1,\ldots,\lambda_j^s
\right)\neq0$.
For any $\bar{u}\in U$ and $t>0$,
\[
H_0+\sum_{\ell=1}^{s}
\left(\bar{u}_\ell+t\lambda_j^\ell\right)H_\ell
=
H_0+\sum_{\ell=1}^{s}\bar{u}_\ell H_\ell
+t\sum_{\ell=1}^{s}\lambda_j^\ell H_\ell
\succeq0.
\]
Hence, $\bar{u}+td_j\in U$ for all $t>0$. 
Since $d_j\neq0$, one has
$\|\bar{u}+td_j\|\to\infty$ as $t\to\infty$, contradicting the
boundedness of $U$. Thus, $\lambda_j^\ell=0,\ \ell=1,\ldots,s$. 
Choosing any $u_2^j\in U$ then gives $\lambda_j^\ell = \lambda_ju_{2,\ell}^j$.

Consequently, for each $r=1,\ldots,p$ and $j=1,\ldots,q$, the vectors $u_1^r,u_2^j\in U$ constructed above satisfy 
\[
\alpha_r^\ell=\alpha_ru_{1,\ell}^r,
\
\lambda_j^\ell=\lambda_ju_{2,\ell}^j,
\
\ell=1,\ldots,s.
\]
By the affine dependence of $f_r$ and $g_j$ on $u$, together with the relations established above, condition \eqref{eq:necessary-polynomial-condition} can be rewritten as
\[
\Phi(x)
=
\sum_{r=1}^{p}\alpha_r f_r(x,u_1^r)
+
\sum_{j=1}^{q}\lambda_j g_j(x,u_2^j)
-
\sum_{r=1}^{p}\alpha_rF_r(x^*,y^{0*},W^*),
\
\Phi(x)\geq0,\ \forall x\in\mathbb{R}^n.
\]
Using the notation $W_\ell$ for the $\ell$-th column of $W$, the balance conditions \eqref{eq:necessary-balance-conditions} yield
\[
\begin{aligned}
&\sum_{r=1}^{p}\alpha_rv_r^\top\left(y^0+Wu_1^r\right)
+
\sum_{j=1}^{q}\lambda_j\theta_j^\top\left(y^0+Wu_2^j\right)\\
&=
\left(
\sum_{r=1}^{p}\alpha_rv_r
+
\sum_{j=1}^{q}\lambda_j\theta_j
\right)^\top y^0 +
\sum_{\ell=1}^{s}
\left(
\sum_{r=1}^{p}\alpha_ru_{1,\ell}^rv_r
+
\sum_{j=1}^{q}\lambda_ju_{2,\ell}^j\theta_j
\right)^\top W_\ell
=0.
\end{aligned}
\]
Consequently, for every
$(x,y^0,W)\in\mathbb{R}^{n+m+ms}$, we have
\begin{equation*}
\sum_{r=1}^{p}\alpha_r\Big[
f_r(x,u_1^r)+v_r^\top(y^0+Wu_1^r)-F_r(x^*,y^{0*},W^*)\Big]+\sum_{j=1}^{q}\lambda_j\Big[
g_j(x,u_2^j)+\theta_j^\top(y^0+Wu_2^j)\Big]=\Phi(x)\geq 0.
\end{equation*}
Let $(x,y^0,W)\in\mathcal{F}$ be arbitrary. Since $u_2^j\in U$, the robust feasibility of $(x,y^0,W)$ implies that
\[
g_j(x,u_2^j)
+
\theta_j^\top\left(y^0+Wu_2^j\right)
\leq0,
\
j=1,\ldots,q.
\]
Together with $\lambda_j\geq0$, this yields
\[
\sum_{r=1}^{p}\alpha_r
\left[
f_r(x,u_1^r)
+
v_r^\top\left(y^0+Wu_1^r\right)
\right]
\geq
\sum_{r=1}^{p}
\alpha_rF_r(x^*,y^{0*},W^*).
\]
Moreover, by the definition of $F_r$ and $\alpha_r>0$, it follows that
\begin{equation}\label{eq:positive-weighted-optimality}
\sum_{r=1}^{p}\alpha_rF_r(x,y^0,W)
\geq
\sum_{r=1}^{p}\alpha_rF_r(x^*,y^{0*},W^*),
\
\forall(x,y^0,W)\in\mathcal{F}.
\end{equation}

We now complete the proof by establishing the Benson proper efficiency of $(x^*,y^{0*},W^*)$. 
Let
\[
\mathcal{K}
:=
\operatorname{cl}\operatorname{cone}
\left(
\mathcal{Y}+\mathbb{R}_+^p
-F(x^*,y^{0*},W^*)
\right).
\]
It follows from \eqref{eq:positive-weighted-optimality} and
$\alpha\in\mathbb{R}_{++}^p$ that
\[
\alpha^\top
\left(
F(x,y^0,W)+z-F(x^*,y^{0*},W^*)
\right)
\geq0,
\
\forall(x,y^0,W)\in\mathcal{F},\
\forall z\in\mathbb{R}_+^p.
\]
Recalling that $\mathcal{Y}=F(\mathcal{F})$, the linearity and continuity of the functional $d\mapsto\alpha^\top d$ yield $\alpha^\top d\geq0,
\
\forall d\in\mathcal{K}$. Suppose, to the contrary, that $0\neq d\in\mathcal{K}\cap(-\mathbb{R}_+^p)$. Since $\alpha\in\mathbb{R}_{++}^p$, we have $\alpha^\top d<0$, which
contradicts the preceding inequality. Hence,
\[
\mathcal{K}\cap(-\mathbb{R}_+^p)=\{0\}.
\]
By Definition~\ref{def:Benson-proper}, $(x^*,y^{0*},W^*)$ is a Benson properly efficient solution of problem \eqref{RTCP}.
The equivalence of Geoffrion and Benson proper efficiency
then implies that it is also Geoffrion properly efficient. This completes the proof.  
\end{proof}

The following example shows that the Slater condition \eqref{eq:slater-condition} in Theorem~\ref{thm:necessary-optimality} cannot be omitted.

\begin{example}
\label{ex:role-slater}
{\rm 

Consider the following uncertain two-stage convex polynomial multiobjective optimization problem:
\begin{align*}
\min\limits_{x\in\mathbb{R}^{2},\,y(\cdot)}\quad \,
& \Bigl(
x_1+x_2^4-u_1,\,
2x_1+u_2
\Bigr)
\tag{\rm EU1}\label{EU1}\\
\text{s.t.}\qquad
& u_1x_1+(u_2+1)x_2\leq0,\ x_1+y(u)\leq2.
\end{align*}
Here, $x=(x_1,x_2)\in\mathbb{R}^2$ is the first-stage decision variable, $y(\cdot): U \to \mathbb{R}$ is the second-stage decision rule, and $u=(u_1,u_2)\in \mathbb{R}^2$ is the uncertain parameter. 
The uncertainty set is given by
\[
U
:=
\left\{
u=(u_1,u_2)^\top\in\mathbb{R}^2
\;\middle|\;
\frac{u_1^2}{4}+u_2^2\leq1
\right\}.
\]
Furthermore, the set $U$ admits the following spectrahedral representation:
\[
U
=
\left\{
u\in\mathbb{R}^2
\;\middle|\;
H_0+u_1H_1+u_2H_2\succeq0
\right\},
\]
where
\[
H_0=
\begin{pmatrix}
4&0&0\\
0&1&0\\
0&0&1
\end{pmatrix},
\quad
H_1=
\begin{pmatrix}
0&0&1\\
0&0&0\\
1&0&0
\end{pmatrix},
\quad
H_2=
\begin{pmatrix}
0&0&0\\
0&0&1\\
0&1&0
\end{pmatrix}.
\]

We use the affine decision rule
$y(u):=y^0+Wu, \ u\in U$, where $y^0\in\mathbb{R}$ and $W\in\mathbb{R}^{1\times2}$ are the decision-rule coefficients. 
Then the corresponding robust counterpart of problem \eqref{EU1} is
\begin{align*}
\min\limits_{\substack{x\in\mathbb{R}^{2},\,y^0\in\mathbb{R},\\ W\in\mathbb{R}^{1\times2}}}
\quad
&\Biggl(
\max_{u\in U}\bigl\{x_1+x_2^4-u_1\bigr\},
\max_{u\in U}\bigl\{2x_1+u_2\bigr\}
\Biggr)
\tag{\rm ER1}\label{ER1}\\
\text{s.t.}\qquad \,
&u_1x_1+(u_2+1)x_2\leq0,
\ \forall u\in U,\\
&x_1+y^0+Wu\leq2,
\ \forall u\in U.
\end{align*}
Problem \eqref{ER1} is a special case of problem \eqref{RTCP}, with $v_1=v_2=\theta_1=0$ and $\theta_2=1$. 
The corresponding polynomial coefficients are
\begin{equation*}
	\begin{alignedat}{3}
f_1^0(x)&=x_1+x_2^4,
&\quad f_1^1(x)&=-1,
&\quad f_1^2(x)&=0,\\
f_2^0(x)&=2x_1,
&\quad f_2^1(x)&=0,
&\quad f_2^2(x)&=1,\\
g_1^0(x)&=x_2,
&\quad g_1^1(x)&=x_1,
&\quad g_1^2(x)&=x_2,\\
g_2^0(x)&=x_1-2,
&\quad g_2^1(x)&=0,
&\quad g_2^2(x)&=0.
\end{alignedat}
\end{equation*}

Consider the feasible point $x^*=(0,0),
\
y^{0*}=0$ and
$W^*=0_{1\times2}$ for problem~\eqref{ER1}. Its robust objective values are
\[
F_1(x^*,y^{0*},W^*)=\max_{u\in U}(-u_1)=2,
\
F_2(x^*,y^{0*},W^*)=\max_{u\in U}u_2=1.
\]
We next show that $(x^*,y^{0*},W^*)$ is a Geoffrion properly efficient solution of problem~\eqref{ER1}. Let $(x,y^0,W)$ be any feasible solution of problem~\eqref{ER1}.
The first robust constraint implies that
\[
\max_{u\in U}
\left\{
u_1x_1+(u_2+1)x_2
\right\}
\leq 0.
\]
Using the support function of $U$, the first robust constraint yields $\sqrt{4x_1^2+x_2^2}+x_2\leq 0$. 
Hence, $x_2\leq0$ and squaring $\sqrt{4x_1^2+x_2^2}\leq -x_2$ yields $x_1=0$.
Consequently, every feasible solution satisfies
\[
F_1(x,y^0,W)
=
x_2^4+2
\geq
2
=
F_1(x^*,y^{0*},W^*),
\
F_2(x,y^0,W)
=
1
=
F_2(x^*,y^{0*},W^*).
\]
Thus, $(x^*,y^{0*},W^*)$ is efficient and no feasible solution can strictly improve either objective. 
Hence, the bounded trade-off condition in definition \ref{def:Geoffrion-proper} is satisfied, and $(x^*,y^{0*},W^*)$ is a Geoffrion properly efficient solution of problem~\eqref{ER1}.

Note that $(0,-1)\in U$. 
For this uncertainty realization, 
\begin{equation*}
g_{1}(x,(0,-1)) =u_{1}x_{1}+\left(u_{2}+1\right) x_{2} =0, \ \forall x \in \mathbb{R}^{2} .
\end{equation*}
Thus, the first robust constraint cannot be strictly satisfied, and the Slater condition \eqref{eq:slater-condition} fails. 

It remains to show conditions \eqref{eq:necessary-lmi-conditions}--\eqref{eq:necessary-polynomial-condition} of
Theorem \ref{thm:necessary-optimality} cannot hold at $(x^*,y^{0*},W^*)$. 
Suppose, to the contrary, that there exist $\alpha=(\alpha_1,\alpha_2)\in\mathbb{R}_{++}^2,
\
\lambda=(\lambda_1,\lambda_2)\in\mathbb{R}_+^2$ and $\alpha_r^\ell,
\
\lambda_j^\ell\in\mathbb{R},
\
r,j,\ell=1,2,$ satisfying these conditions.
The balance conditions in \eqref{eq:necessary-balance-conditions}, together with $v_1=v_2=\theta_1=0$ and $\theta_2=1$, give
\begin{equation}
\lambda_2=0,
\
\lambda_2^1=\lambda_2^2=0.
\label{eq:example-balance}
\end{equation}
The matrix inequality in condition \eqref{eq:necessary-lmi-conditions} associated with $\lambda_1$ is $\lambda_1H_0+\lambda_1^1H_1+\lambda_1^2H_2\succeq0$, or equivalently,
\begin{equation}
\begin{pmatrix}
4\lambda_1 & 0 & \lambda_1^1\\
0 & \lambda_1 & \lambda_1^2\\
\lambda_1^1 & \lambda_1^2 & \lambda_1
\end{pmatrix}
\succeq0.
\label{eq:example-lmi}
\end{equation}
If $\lambda_1>0$, applying the Schur complement to
\eqref{eq:example-lmi} yields
$\frac{(\lambda_1^1)^2}{4}+
(\lambda_1^2)^2\leq(\lambda_1)^2.$
If $\lambda_1=0$, the positive semidefiniteness in \eqref{eq:example-lmi} directly implies $\lambda_1^1=\lambda_1^2=0$. Thus, in either case,
\begin{equation}
\frac{(\lambda_1^1)^2}{4}
+
(\lambda_1^2)^2
\leq
(\lambda_1)^2.
\label{eq:example-lmi-inequality}
\end{equation}
Substituting the polynomial coefficients and \eqref{eq:example-balance} into \eqref{eq:necessary-polynomial-condition}, we obtain
\begin{equation}
\begin{aligned}
\Phi(x)
=
\alpha_1x_2^4
+
\left(
\alpha_1+2\alpha_2+\lambda_1^1
\right)x_1+
\left(
\lambda_1+\lambda_1^2
\right)x_2
-\alpha_1^1+\alpha_2^2-2\alpha_1-\alpha_2.
\end{aligned}
\label{eq:example-polynomial}
\end{equation}
By condition \eqref{eq:necessary-polynomial-condition}, $x^{*}=(0,0)$ is a global minimizer of $\Phi$. Hence, $\nabla \Phi\left(x^{*}\right)=0$. Together with \eqref{eq:example-polynomial}, this yields
\begin{equation*}
\alpha_{1}+2 \alpha_{2}+\lambda_{1}^{1}=0, \quad \lambda_{1}+\lambda_{1}^{2}=0 .
\end{equation*}
Substituting  $\lambda_{1}^{2}=-\lambda_{1}$ into \eqref{eq:example-lmi-inequality}, we obtain
\begin{equation*}
\frac{\left(\lambda_{1}^{1}\right)^{2}}{4}+\left(\lambda_{1}\right)^{2} \leq\left(\lambda_{1}\right)^{2} .
\end{equation*}
Thus, $\lambda_{1}^{1}=0$, and $\alpha_{1}+2 \alpha_{2}=0$, which contradicts $\alpha_1>0$ and $\alpha_2>0$.
This shows that the necessity result may fail without the Slater condition.
}\end{example}

The optimality conditions established in this section apply to general convex polynomial data. 
In the next section, we consider SOS-convex polynomial data and show that the global polynomial nonnegativity condition~\eqref{eq:necessary-polynomial-condition} admits an exact SOS representation.
We then develop the corresponding moment-SDP formulations.

\section{SOS Characterizations and Moment-SDP Formulations}\label{sect:4}
We first recall some notation about the polynomial.
A polynomial $h\in\mathbb{R}[x]$ is called a sum-of-squares (SOS) polynomial if there exist polynomials
$h_{1},\dots,h_{N}\in\mathbb{R}[x]$, such that
\[
h(x)=\sum_{k=1}^{N}h_k(x)^2.
\]
The cone of all SOS polynomials in $\mathbb{R}[x]$ is denoted by
$\Sigma[x]$.

\begin{definition}[SOS-convex polynomial]\label{def:SOS-convex}
A real polynomial $h:\mathbb{R}^n\to\mathbb{R}$ is called SOS-convex if the polynomial
\[
H(x,z):=h(x)-h(z)-\nabla h(z)^\top(x-z)
\]
is an SOS polynomial in variables $(x,z)$.
\end{definition}

By \cite[Theorem 3.1]{ahmadi2013complete}, this definition is equivalent to
requiring that $\nabla^2 h(x)$ be an SOS matrix, that is, $z^\top\nabla^2 h(x)z \in \Sigma[x,z]$. 
The following lemma provides a sufficient condition for an SOS-convex polynomial to be SOS. 

\begin{lemma}\label{lem:SOS-convex}
Let $h:\mathbb{R}^n\to\mathbb{R}$ be an SOS-convex polynomial. Suppose that there exists $\bar{x}\in\mathbb{R}^n$ such that $h(\bar{x})=0$ and
$\nabla h(\bar{x})=0$.
Then $h\in\Sigma[x]$.
\end{lemma}

We consider SOS-convex polynomial data. 
Under this additional assumption, condition~\eqref{eq:necessary-polynomial-condition} admits an exact SOS characterization, as shown in the following theorem.

\begin{theorem}[SOS characterization under SOS-convexity]
\label{thm:SOS-characterization}
Let $(x^*,y^{0*},W^*)\in\mathcal{F}$. 
Suppose that $f_r(\cdot,u)$, $r=1,\ldots,p$, and $g_j(\cdot,u)$,
$j=1,\ldots,q$, are SOS-convex for every $u\in U$, and that the Slater condition \eqref{eq:slater-condition} holds. 
Then $(x^*,y^{0*},W^*)$ is a Geoffrion properly efficient solution of problem~\eqref{RTCP} if and only if there exist $\alpha\in\mathbb{R}_{++}^p$,
$\lambda\in\mathbb{R}_+^q$, and
$\alpha_r^\ell \ \lambda_j^\ell\in\mathbb{R}$, for $r=1,\ldots,p$,
$j=1,\ldots,q$, $\ell=1,\ldots,s$, satisfying
conditions \eqref{eq:necessary-lmi-conditions} and \eqref{eq:necessary-balance-conditions}, such that the polynomial $\Phi$ defined in  \eqref{eq:necessary-polynomial-condition} satisfies
\[
\Phi\in\Sigma[x],
\
\Phi(x^*)=0.
\]
\end{theorem}

\begin{proof}
For sufficiency, $\Phi\in\Sigma[x]$ and $\Phi(x^*)=0$ implies
condition~\eqref{eq:necessary-polynomial-condition}. 
The conclusion follows from Theorem~\ref{thm:sufficient-optimality}. 

For the necessity implication, suppose that
$(x^*,y^{0*},W^*)$ is a Geoffrion properly efficient solution of
problem~\eqref{RTCP}. By Theorem~\ref{thm:necessary-optimality}, there exist multipliers satisfying conditions~\eqref{eq:necessary-lmi-conditions}--\eqref{eq:necessary-polynomial-condition}. The construction in the proof of Theorem~\ref{thm:sufficient-optimality} yields vectors $u_1^r\in U$, $r=1,\ldots,p$, and $u_2^j\in U$, $j=1,\ldots,q$, such that
\[
\alpha_r^\ell=\alpha_r u_{1,\ell}^r,
\
\lambda_j^\ell=\lambda_j u_{2,\ell}^j,
\
\ell=1,\ldots,s.
\]
Using the affine dependence of $f_r$ and $g_j$ on $u$, we can write $\Phi$ as
\[
\Phi(x)
=
\sum_{r=1}^p \alpha_r f_r(x,u_1^r)
+
\sum_{j=1}^q \lambda_j g_j(x,u_2^j)
-
\sum_{r=1}^p
\alpha_r F_r(x^*,y^{0*},W^*).
\]
By assumption, $f_r(\cdot,u_1^r)$ and $g_j(\cdot,u_2^j)$ are
SOS-convex. 
Since $\alpha_r>0$ and $\lambda_j\geq0$, the above representation shows that $\Phi$ is SOS-convex.
Condition~\eqref{eq:necessary-polynomial-condition} implies that $x^*$ is a global minimizer of $\Phi$ with $\Phi(x^*)=0$. Hence, $\nabla\Phi(x^*)=0$ and lemma~\ref{lem:SOS-convex} yields $\Phi\in\Sigma[x]$.
\end{proof}

It is worth noting that the Slater condition \eqref{eq:slater-condition} is required only for necessity, not for sufficiency.
We next develop a moment-SDP approach for computing Geoffrion properly efficient solutions via strictly positive weighted-sum scalarization.

For each vector
$\alpha=(\alpha_1,\ldots,\alpha_p)\in\mathbb{R}_{++}^p$,
consider the following weighted-sum scalarization of problem~\eqref{RTCP}:
\begin{align*}
\min\limits_{\substack{x\in\mathbb{R}^{n},\,y^0\in\mathbb{R}^{m},\\ W\in\mathbb{R}^{m\times s}}}
\quad \
&
\sum_{r=1}^{p}\alpha_r
\max_{u\in U}
\Bigl\{
f_r(x,u)+v_r^\top\bigl(y^0+Wu\bigr)
\Bigr\}
\tag{$\rm RTCP_{\alpha}$}\label{RTCPalpha}
\\
\text{s.t.}\quad \qquad
&
g_j(x,u)+\theta_j^\top\bigl(y^0+Wu\bigr)\leq 0,
\
\forall u\in U,\ j=1,\ldots,q.
\end{align*}

The following proposition establishes the equivalence between Geoffrion properly efficiency and optimality for a weighted-sum scalarization with strictly positive weights.

\begin{proposition}\label{prop:weighted-sum}
Let $(x^*,y^{0*},W^*)\in\mathcal{F}$. Then
$(x^*,y^{0*},W^*)$ is a Geoffrion properly efficient solution of
problem~\eqref{RTCP} if and only if there exists
$\alpha\in\mathbb{R}_{++}^p$ such that $(x^*,y^{0*},W^*)$ is an optimal solution of problem~\eqref{RTCPalpha}.
\end{proposition}

\begin{proof}
The necessity follows from
the equivalence of Geoffrion and Benson proper efficiency and Lemma~\ref{lem:positive-supporting-vector}. 
It remains only to prove the sufficiency. Suppose that $(x^*,y^{0*},W^*)$ is an optimal solution of problem~\eqref{RTCPalpha} for some $\alpha\in\mathbb{R}_{++}^p$.
Then
\[
\sum_{r=1}^{p}\alpha_r F_r(x,y^0,W)
\geq
\sum_{r=1}^{p}\alpha_r F_r(x^*,y^{0*},W^*),
\
\forall (x,y^0,W)\in\mathcal{F}.
\]
Following the same argument as in the proof of
Theorem~\ref{thm:sufficient-optimality}, $(x^*,y^{0*},W^*)$ is a Benson properly efficient and hence Geoffrion properly efficient.
\end{proof}

By proposition~\ref{prop:weighted-sum}, Geoffrion properly efficient solutions of problem~\eqref{RTCP} can be obtained by  solving the strictly positive weighted-sum problem~\eqref{RTCPalpha}. 
Assume that $f_r(\cdot,u)$ and $g_j(\cdot,u)$ are SOS-convex for every $u\in U$. 
We now construct a finite-dimensional moment-SDP relaxation of problem~\eqref{RTCPalpha}.

Choose $d\in\mathbb{N}_{+}$ sufficiently large so that
$f_r^\ell,g_j^\ell\in\mathbb{R}[x]_{2d}$ for all
$r=1,\ldots,p$, $j=1,\ldots,q$, and $\ell=0,\ldots,s$,
where $\mathbb{R}[x]_{2d}$ denotes the space of real polynomials in $x$ of degree at most $2d$.
For a multi-index
$\beta=(\beta_1,\ldots,\beta_n)\in\mathbb{N}^n$, where
$\mathbb{N}:=\{0,1,\ldots\}$, write
\[
|\beta|:=\sum_{i=1}^{n}\beta_i,
\
\mathbb{N}_d^n
:=
\left\{
\beta\in\mathbb{N}^n:|\beta|\leq d
\right\},
\
x^\beta:=x_1^{\beta_1}\cdots x_n^{\beta_n}.
\]
Denote by $z_d(x):=
\bigl(x^\beta\bigr)_{\beta\in\mathbb{N}_d^n}$ the column vector of all monomials in $x$ of degree at most $d$. 
Let $\mathbf{m}=\bigl(m_\beta\bigr)_{\beta\in\mathbb{N}_{2d}^n}$ be a truncated moment sequence. 
Its associated moment matrix $M_d(\mathbf{m})$ is defined by
\[
\bigl[M_d(\mathbf{m})\bigr]_{\beta,\gamma}
:=
m_{\beta+\gamma},
\
\beta,\gamma\in\mathbb{N}_d^n.
\]
For $h(x)=
\sum_{\beta\in\mathbb{N}_{2d}^n}h_\beta x^\beta
\in\mathbb{R}[x]_{2d},$ define 
the associated Riesz linear functional by
\[
L_{\mathbf{m}}(h)
:=
\sum_{\beta\in\mathbb{N}_{2d}^n}
h_\beta m_\beta.
\]

We assume that $U$ satisfies the following strict feasibility condition: there exists $\widehat{u}\in\mathbb{R}^{s}$ such that
\begin{equation}\label{eq:strict-feasibility-U}
H_0+\sum_{\ell=1}^{s}\widehat{u}_{\ell}H_{\ell}\succ0.
\end{equation}
Under condition~\eqref{eq:strict-feasibility-U}, SDP strong duality (see~\cite{vandenberghe1996semidefinite}) yields the following exact representation. 
For any $(a_0,a_1,\ldots,a_s)\in\mathbb{R}^{s+1}$,
\begin{equation}\label{eq:spectrahedral-support}
\begin{aligned}
\sup_{u\in U}
\left\{
a_0+\sum_{\ell=1}^{s}a_\ell u_\ell
\right\}
=
\min_{Z\succeq0}
\left\{
a_0+\langle H_0,Z\rangle
\ \middle|\
a_\ell+\langle H_\ell,Z\rangle=0,\ 
\ell=1,\ldots,s
\right\}.
\end{aligned}
\end{equation}

Let $W_\ell\in\mathbb{R}^m$ denote the $\ell$-th column of $W$,
$\ell=1,\ldots,s$. 
Replacing the polynomial terms by their associated Riesz linear functionals and applying~\eqref{eq:spectrahedral-support} separately
to the robust objective and constraints, we formulate the following moment-SDP relaxation of problem~\eqref{RTCPalpha}.
\begin{align*}
\inf\limits_{\substack{\mathbf{m},\,y^0,\,W,\\ Z_r,\,S_j}}\quad \
&
\sum_{r=1}^{p}\alpha_r
\Bigl[
L_{\mathbf{m}}(f_r^0)
+v_r^{\top}y^0
+\langle H_0,Z_r\rangle
\Bigr]
\tag{$\rm RTCP_{\alpha}^{*}$}\label{RTCPalphastar}\\
\text{s.t.}\qquad
&
L_{\mathbf{m}}(f_r^\ell)
+v_r^{\top}W_\ell
+\langle H_\ell,Z_r\rangle
=0,
\ r=1,\ldots,p,\ \ell=1,\ldots,s,
\\
&
L_{\mathbf{m}}(g_j^0)
+\theta_j^{\top}y^0
+\langle H_0,S_j\rangle
\leq 0,
\ j=1,\ldots,q,
\\
&
L_{\mathbf{m}}(g_j^\ell)
+\theta_j^{\top}W_\ell
+\langle H_\ell,S_j\rangle
=0,
\ j=1,\ldots,q,\ \ell=1,\ldots,s,
\\
&
m_{\mathbf{0}}=1,\
M_d(\mathbf{m})\succeq 0,
\\
&
Z_r\succeq 0,\ r=1,\ldots,p,
\\
&
S_j\succeq 0,\ j=1,\ldots,q.
\end{align*}
Here, $\mathbf{m}\in\mathbb{R}^{|\mathbb{N}_{2d}^{n}|}$,
$y^0\in\mathbb{R}^{m}$, $W\in\mathbb{R}^{m\times s}$,
$Z_r\in\mathbb{S}^{k}$, $r=1,\ldots,p$, and
$S_j\in\mathbb{S}^{k}$, $j=1,\ldots,q$.

The following Jensen-type inequality for SOS-convex polynomials will be used to recover a first-stage decision from the moment variables.

\begin{lemma}[\cite{lasserre2009convexity,liguo2021finding}]
\label{lem:sos-convex-jensen}
Let $h\in\mathbb{R}[x]_{2d}$ be SOS-convex, and let  $\mathbf{m}=\left(m_{\beta}\right)_{\beta \in \mathbb{N}_{2d}^{n}}$ satisfy  $m_{\mathbf{0}}=1$ and $M_d(\mathbf{m})\succeq0$.  Then
\[
h\bigl(x(\mathbf{m})\bigr)
\leq
L_{\mathbf{m}}(h),
\]
where $x(\mathbf{m}):=\left(L_{\mathbf{m}}\left(x_{1}\right), \ldots, L_{\mathbf{m}}\left(x_{n}\right)\right).$
\end{lemma}

The following theorem
establishes an exact moment-SDP relaxation and provides a correspondence between its optimal solutions and the Geoffrion properly efficient solutions of problem~\eqref{RTCP}.

\begin{theorem}\label{thm:moment-sdp}
Suppose that $f_r(\cdot,u)$, $r=1,\ldots,p$, and $g_j(\cdot,u)$, $j=1,\ldots,q$, are SOS-convex for every $u\in U$, and that the spectrahedral representation of $U$ satisfies the strict feasibility condition~\eqref{eq:strict-feasibility-U}. Then the following statements hold.

\medskip
\noindent
\textnormal{\rm(i)}
For every $\alpha\in\mathbb{R}_{++}^{p}$,
\[
\operatorname{val}(\mathrm{RTCP}_{\alpha})
=
\operatorname{val}(\mathrm{RTCP}_{\alpha}^{*}).
\]

\medskip
\noindent
\textnormal{\rm(ii)}
Let $(x^*,y^{0*},W^*)$ be a Geoffrion properly efficient solution of problem~\eqref{RTCP}. Then there exist
$\alpha\in\mathbb{R}_{++}^{p}$, $(Z_r^*)_{r=1}^{p}$, and
$(S_j^*)_{j=1}^{q}$ such that
\[
\left(
\mathbf{m}^{x^*},y^{0*},W^*,
(Z_r^*)_{r=1}^{p},(S_j^*)_{j=1}^{q}
\right)
\]
is an optimal solution of problem~\eqref{RTCPalphastar}, where
\[
\mathbf{m}^{x^*}
:=
\left(
(x^*)^\beta
\right)_{\beta\in\mathbb{N}_{2d}^{n}}.
\]

\medskip
\noindent
\textnormal{\rm(iii)}
Fix $\alpha\in\mathbb{R}_{++}^{p}$ and let
\[
\left(
\mathbf{m}^*,y^{0*},W^*,
(Z_r^*)_{r=1}^{p},(S_j^*)_{j=1}^{q}
\right)
\]
be an optimal solution of problem~\eqref{RTCPalphastar}. Then
$(\bar{x},y^{0*},W^*)$ is a Geoffrion properly efficient
solution of problem~\eqref{RTCP}, where
\[
\bar{x}
:=
x(\mathbf{m}^*)
=
\left(
L_{\mathbf{m}^*}(x_1),\ldots,
L_{\mathbf{m}^*}(x_n)
\right).
\]
\end{theorem}

\begin{proof}
\textnormal{\rm(i)}
We first show that $\operatorname{val}(\mathrm{RTCP}_{\alpha}^{*})
\leq
\operatorname{val}(\mathrm{RTCP}_{\alpha}).$
Let $(x,y^0,W)$ be any feasible solution of problem~\eqref{RTCPalpha} and
define $\mathbf{m}^{x}:=(x^\beta)_{\beta\in\mathbb{N}_{2d}^{n}}$.
Then $m_{\mathbf{0}}^{x}=1,
\
M_d(\mathbf{m}^{x})
=
z_d(x)z_d(x)^\top
\succeq0$.
For each $r=1,\ldots,p$, applying
\eqref{eq:spectrahedral-support} to the affine function
\[
u\longmapsto
f_r^0(x)+v_r^\top y^0
+
\sum_{\ell=1}^{s}u_\ell
\left(
f_r^\ell(x)+v_r^\top W_\ell
\right)
\]
yields a matrix $Z_r\succeq0$ such that
\[
L_{\mathbf{m}^{x}}(f_r^\ell)
+v_r^{\top}W_\ell
+\langle H_\ell,Z_r\rangle
=0,
\
\ell=1,\ldots,s,
\]
\[
F_r(x,y^0,W)
=
L_{\mathbf{m}^{x}}(f_r^0)
+v_r^{\top}y^0
+\langle H_0,Z_r\rangle.
\]
Similarly, for each $j=1,\ldots,q$ there exists $S_j\succeq0$ such that
\[
L_{\mathbf{m}^{x}}(g_j^\ell)
+\theta_j^{\top}W_\ell
+\langle H_\ell,S_j\rangle
=0,
\
\ell=1,\ldots,s,
\]
\[
L_{\mathbf{m}^{x}}(g_j^0)
+\theta_j^{\top}y^0
+\langle H_0,S_j\rangle
\leq0.
\]
Therefore, $\bigl(
\mathbf{m}^{x},y^0,W,
(Z_r)_{r=1}^{p},
(S_j)_{j=1}^{q}
\bigr)$ is feasible for problem~\eqref{RTCPalphastar}, with objective value 
\[
\begin{aligned}
&
\sum_{r=1}^{p}\alpha_r
\left[
L_{\mathbf{m}^{x}}(f_r^0)
+v_r^{\top}y^0
+\langle H_0,Z_r\rangle
\right]=
\sum_{r=1}^{p}\alpha_rF_r(x,y^0,W).
\end{aligned}
\]
Since $\left(x, y^{0}, W\right)$ was arbitrary, taking the infimum gives $\operatorname{val}(\mathrm{RTCP}_{\alpha}^{*})
\leq
\operatorname{val}(\mathrm{RTCP}_{\alpha})$.

Conversely, let $\bigl(
\mathbf{m},y^0,W,(Z_r)_{r=1}^{p},(S_j)_{j=1}^{q}
\bigr)$
be any feasible solution of problem~\eqref{RTCPalphastar}, and let $\bar{x}:=x(\mathbf{m})$.
For any $u\in U$ and $j=1,\ldots,q$, Lemma~\ref{lem:sos-convex-jensen}
and the constraints of problem~\eqref{RTCPalphastar} give
\[
\begin{aligned}
g_j(\bar{x},u)+\theta_j^\top(y^0+Wu) &\leq \,
L_{\mathbf{m}}\bigl(g_j(\cdot,u)\bigr)
+\theta_j^\top(y^0+Wu)
\\
&=\,
L_{\mathbf{m}}(g_j^0)+\theta_j^\top y^0
+\sum_{\ell=1}^{s}u_\ell
\left[
L_{\mathbf{m}}(g_j^\ell)
+\theta_j^\top W_\ell
\right]
\\
&=\,
L_{\mathbf{m}}(g_j^0)+\theta_j^\top y^0
+\langle H_0,S_j\rangle
-
\left\langle
H_0+\sum_{\ell=1}^{s}u_\ell H_\ell,
S_j
\right\rangle
\\
& \leq\, 0,
\end{aligned}
\]
where the last inequality follows from
\begin{equation*}
H_{0}+\sum_{\ell=1}^{s} u_{\ell} H_{\ell} \succeq 0, \quad S_{j} \succeq 0, \quad L_{\mathbf{m}}\left(g_{j}^{0}\right)+\theta_{j}^{\top} y^{0}+\left\langle H_{0}, S_{j}\right\rangle \leq 0 .
\end{equation*}
Hence, $(\bar{x},y^0,W)$ is feasible for problem~\eqref{RTCPalpha}.

Similarly, for any $u\in U$ and $r=1,\ldots,p$,
\[
\begin{aligned}
f_r(\bar{x},u)+v_r^\top(y^0+Wu)
&\leq \,
L_{\mathbf{m}}\bigl(f_r(\cdot,u)\bigr)
+v_r^\top(y^0+Wu)
\\
&=\,
L_{\mathbf{m}}(f_r^0)+v_r^\top y^0
+\langle H_0,Z_r\rangle
-
\left\langle
H_0+\sum_{\ell=1}^{s}u_\ell H_\ell,
Z_r
\right\rangle
\\
&\leq\,
L_{\mathbf{m}}(f_r^0)+v_r^\top y^0
+\langle H_0,Z_r\rangle.
\end{aligned}
\]
Taking the supremum over $u\in U$ gives
\begin{equation}\label{eq:objective-recovery}
F_r(\bar{x},y^0,W)
\leq
L_{\mathbf{m}}(f_r^0)
+v_r^\top y^0
+\langle H_0,Z_r\rangle,
\ r=1,\ldots,p.
\end{equation}
Therefore,
\[
\sum_{r=1}^{p}\alpha_rF_r(\bar{x},y^0,W)
\leq
\sum_{r=1}^{p}\alpha_r
\left[
L_{\mathbf{m}}(f_r^0)+v_r^\top y^0
+\langle H_0,Z_r\rangle
\right],
\]
it follows that $\operatorname{val}(\mathrm{RTCP}_{\alpha})
\leq
\operatorname{val}(\mathrm{RTCP}_{\alpha}^{*})$. This proves assertion~\textnormal{\rm(i)}.
 
\noindent
\textnormal{\rm(ii)}
By Proposition~\ref{prop:weighted-sum}, there exists
$\alpha\in\mathbb{R}_{++}^{p}$ such that $(x^*,y^{0*},W^*)$ is an optimal solution of problem~\eqref{RTCPalpha}. By the lifting construction in the proof of~\textnormal{\rm(i)}, there exist
$Z_r^*\succeq0$, $r=1,\ldots,p$, and
$S_j^*\succeq0$, $j=1,\ldots,q$, such that $\left(
\mathbf{m}^{x^*},y^{0*},W^*,
(Z_r^*)_{r=1}^{p},(S_j^*)_{j=1}^{q}
\right)$
is feasible for problem~\eqref{RTCPalphastar}, with objective value
\[
\sum_{r=1}^{p}\alpha_r
\left[
L_{\mathbf{m}^{x^*}}(f_r^0)
+v_r^\top y^{0*}
+\langle H_0,Z_r^*\rangle
\right]=
\sum_{r=1}^{p}\alpha_r
F_r(x^*,y^{0*},W^*)
=
\operatorname{val}(\mathrm{RTCP}_{\alpha}),
\]
where the last equality follows from {\rm(i)}.
Hence the constructed solution is optimal for problem~\eqref{RTCPalphastar}.

\noindent
\textnormal{\rm(iii)}
By the recovery argument in the proof of ~\textnormal{\rm(i)}, $(\bar{x},y^{0*},W^*)$ is feasible for problem~\eqref{RTCPalpha}. Applying \eqref{eq:objective-recovery} to the given optimal solution of problem~\eqref{RTCPalphastar}, we obtain
\[
F_r(\bar{x},y^{0*},W^*)
\leq
L_{\mathbf{m}^*}(f_r^0)
+v_r^\top y^{0*}
+\langle H_0,Z_r^*\rangle,
\ r=1,\ldots,p.
\]
Consequently,
\[
\begin{aligned}
\operatorname{val}(\mathrm{RTCP}_{\alpha}) &\,\leq
\sum_{r=1}^{p}\alpha_r
F_r(\bar{x},y^{0*},W^*)
\leq\,
\sum_{r=1}^{p}\alpha_r
\left[
L_{\mathbf{m}^*}(f_r^0)
+v_r^\top y^{0*}
+\langle H_0,Z_r^*\rangle
\right]
\\
&=\,
\operatorname{val}(\mathrm{RTCP}_{\alpha}^{*})=\,
\operatorname{val}(\mathrm{RTCP}_{\alpha}).
\end{aligned}
\]
Hence, $(\bar{x},y^{0*},W^*)$ is an optimal
solution of problem~\eqref{RTCPalpha}. Since
$\alpha\in\mathbb{R}_{++}^{p}$, Proposition~\ref{prop:weighted-sum} implies that $(\bar{x},y^{0*},W^*)$ is a Geoffrion properly efficient solution of problem~\eqref{RTCP}.
\end{proof}

We now present a numerical example to illustrate how the moment-SDP formulation developed above can be used to compute Geoffrion properly efficient solutions of a robust two-stage convex polynomial multiobjective optimization problem.

\begin{example}
{\rm
Consider the following adjustable robust two-stage multiobjective polynomial optimization problem:
\begin{align*}
\min\limits_{\substack{x\in\mathbb{R}_+^3,\,y^0\in\mathbb{R},\\
W\in\mathbb{R}^{1\times 2}}}
\quad
&
\left(
\max_{u\in U}
\left\{
f_1(x,u)+v_1^\top y(u)
\right\},
\max_{u\in U}
\left\{
f_2(x,u)+v_2^\top y(u)
\right\}
\right)
\tag{$\rm ER2$}\label{ER2}
\\
\text{\rm s.t.}\qquad \,\,
&
g_j(x,u)+\theta_j^\top y(u)\leq 0,
\
j=1,2,3,\
\forall u\in U,
\\
&
y(u)=y^0+Wu,
\
\forall u\in U,
\end{align*}
where
\begin{align*}
f_1(x,u)
={}&
(x_1+x_2+x_3-4)^4+(x_3-2)^2
-u_1x_1-\frac{3}{2}x_2-x_3,
\\
f_2(x,u)
={}&
(x_1+x_2+x_3-4)^4+(x_3-2)^2
-\frac{4}{3}x_1-u_2x_2-\frac{2}{3}x_3,
\\
g_1(x,u)
={}&
u_1x_1+x_2+\frac{1}{2}x_3-4,\
g_2(x,u)
={}
\frac{1}{2}x_1+u_2x_2+x_3-3,
\\
g_3(x,u)
={}&
x_1+x_2+x_3-4.
\end{align*}
The uncertainty set is given by
\[
U
:=
\left\{
u=(u_1,u_2)^\top\in\mathbb{R}^2
\;\middle|\;
(u_1-1)^2+(u_2-1)^2\leq\frac{1}{4}
\right\}.
\]
Furthermore, the uncertainty set $U$ admits the following
spectrahedral representation:
\[
U
=
\left\{
u\in\mathbb{R}^2
\;\middle|\;
H_0+u_1H_1+u_2H_2\succeq0
\right\},
\]
where
\[
H_0
=
\begin{pmatrix}
\frac{1}{2} & -1 & -1\\
-1 & \frac{1}{2} & 0\\
-1 & 0 & \frac{1}{2}
\end{pmatrix},
\
H_1
=
\begin{pmatrix}
0 & 1 & 0\\
1 & 0 & 0\\
0 & 0 & 0
\end{pmatrix},
\
H_2
=
\begin{pmatrix}
0 & 0 & 1\\
0 & 0 & 0\\
1 & 0 & 0
\end{pmatrix}.
\]
%
Taking $\widehat{u}=(1,1)^\top\in U$, we obtain $H_0+\widehat{u}_1H_1+\widehat{u}_2H_2=\frac{1}{2}I_3\succ0.$
Thus, the strict feasible condition \eqref{eq:strict-feasibility-U} is satisfied.
Moreover, for every $u\in U$, the two objective polynomials have the
same Hessian matrix
\[
\nabla_x^2 f_r(x,u)
=
12(x_1+x_2+x_3-4)^2
\begin{pmatrix}
1\\1\\1
\end{pmatrix}
\begin{pmatrix}
1\\1\\1
\end{pmatrix}^{\top}
+
2
\begin{pmatrix}
0\\0\\1
\end{pmatrix}
\begin{pmatrix}
0\\0\\1
\end{pmatrix}^{\top},
\, r=1,2,
\]
which implies $f_1(\cdot,u)$ and $f_2(\cdot,u)$ are SOS-convex polynomials. Since $g_j(\cdot,u)$, $j=1,2,3$, are affine in $x$ for every $u\in U$, they are also SOS-convex.

\noindent\textbf{Case 1 (A two-stage program).} 
Let $v=
\left(-\frac{1}{2},-\frac{1}{2}\right)^\top,
\
\theta
=
(1,1,1)^\top$.
For each weight vector
$\alpha=(\alpha_1,\alpha_2)\in\mathbb{R}_{++}^2$, consider the
following robust two-stage weighted-sum problem:
\begin{align*}
\min\limits_{\substack{x\in\mathbb{R}_+^3,\,y^0\in\mathbb{R},\\ W\in\mathbb{R}^{1\times 2}}}
\quad
&
\alpha_1
\max_{u\in U}
\left\{
f_1(x,u)-\frac{1}{2}y(u)
\right\}
+
\alpha_2
\max_{u\in U}
\left\{
f_2(x,u)-\frac{1}{2}y(u)
\right\}
\tag{$\rm ER2_{\alpha}$}\label{ER2alpha}
\\
\text{\rm s.t.}\quad \quad \ 
&
g_j(x,u)+y(u)\leq0,
\
j=1,2,3,\
\forall u\in U,
\\
&
y(u)=y^0+Wu,
\
\forall u\in U.
\end{align*}
Since the objective polynomials have degree four, we take $d=2$. 
The corresponding moment-SDP formulation of problem  \eqref{ER2alpha} is given by 
%
\begin{align*}
\inf\limits_{\substack{\mathbf{m},\,y^0,\,W, \\Z_r,\,S_j}}
\quad \
&
\sum_{r=1}^{2}\alpha_r
\left[
L_{\mathbf{m}}(f_r^0)
-\frac{1}{2}y^0
+\langle H_0,Z_r\rangle
\right]
\tag{$\rm ER2_{\alpha}^{*}$}\label{ER2alphaStar}
\\
\text{\rm s.t.}\qquad
&
L_{\mathbf{m}}(f_r^\ell)
-\frac{1}{2}W_\ell
+\langle H_\ell,Z_r\rangle
=0,
\
r=1,2,\ \ell=1,2,
\\
&
L_{\mathbf{m}}(g_j^0)
+y^0
+\langle H_0,S_j\rangle
\leq0,
\
j=1,2,3,
\\
&
L_{\mathbf{m}}(g_j^\ell)
+W_\ell
+\langle H_\ell,S_j\rangle
=0,
\
j=1,2,3,\ \ell=1,2,
\\
&
m_{\mathbf{0}}=1,
\
M_2(\mathbf{m})\succeq0,
\\
&
L_{\mathbf{m}}(x_i)\geq0,
\
i=1,2,3,
\\
&
Z_r\succeq0,
\
r=1,2,
\\
&
S_j\succeq0,
\
j=1,2,3,
\end{align*}
where $\mathbf{m}=(m_\beta)_{\beta\in\mathbb{N}_4^3},
\
Z_r\in\mathbb{S}_+^3,\ r=1,2,
\
S_j\in\mathbb{S}_+^3,\ j=1,2,3$.

We choose nine strictly positive weight vectors of the form
\[
\alpha^{(k)}
:=
\left(
\frac{k}{10},
1-\frac{k}{10}
\right),
\
k=1,\ldots,9.
\]
For each $\alpha^{(k)}$,
we solve problem~\eqref{ER2alphaStar} numerically in MATLAB using CVX~\cite{grant2008graph,grant2014cvx}.
Let 
\[
\left(
\mathbf{m}^{*(k)},
y^{0*(k)},
W^{*(k)},
(Z_r^{*(k)})_{r=1}^{2},
(S_j^{*(k)})_{j=1}^{3}
\right)
\]
denote an optimal solution obtained in this way. By Theorem~\ref{thm:moment-sdp}, $\left(
\bar{x}^{(k)},
y^{0*(k)},
W^{*(k)}
\right)$ is a Geoffrion properly efficient solution of problem~\eqref{ER2}, where
\[
\bar{x}^{(k)}
:=
x\bigl(\mathbf{m}^{*(k)}\bigr)
=
\left(
L_{\mathbf{m}^{*(k)}}(x_1),
L_{\mathbf{m}^{*(k)}}(x_2),
L_{\mathbf{m}^{*(k)}}(x_3)
\right),
\]
the associated affine second-stage decision rule is
\[
y^{*(k)}(u)
=
y^{0*(k)}+W^{*(k)}u,
\
\forall u\in U.
\]

We next evaluate the accuracy and feasibility of the recovered numerical solutions. 
Let $v_{\mathrm{SDP}}^{(k)}$ denote the optimal value of problem \eqref{ER2alphaStar} corresponding to $\alpha^{(k)}$. 
The absolute gap between this value and the weighted objective value evaluated at the recovered solution is defined by
\[
\operatorname{gap}^{(k)}
:=
\left|
v_{\mathrm{SDP}}^{(k)}
-
\sum_{r=1}^{2}
\alpha_r^{(k)}
F_r\bigl(\bar{x}^{(k)},y^{0*(k)},W^{*(k)}\bigr)
\right|.
\]
The maximum robust constraint violation is defined by
\[
\operatorname{viol}^{(k)}
:=
\max\left\{
0,\,
\max_{j=1,2,3}
\sup_{u\in U}
\left[
g_j\bigl(\bar{x}^{(k)},u\bigr)
+y^{0*(k)}+W^{*(k)}u
\right]
\right\}.
\]
The recovered first-stage solutions are reported in
Table~\ref{tab:numerical-results}, while the SDP optimal values,
the recovered weighted objective values, the absolute gaps, and the maximum robust constraint violations are reported in Table ~\ref{tab:case1-results}.

\noindent\textbf{Case 2 (A single-stage program).}
For comparison, we consider the single-stage model obtained from Case 1 by fixing $y(u) \equiv 0$.
Equivalently, we set  $v=(0,0)^{\top}$ and  $\theta=(0,0,0)^{\top}$ in problem (ER2). 
For each
$\alpha=(\alpha_1,\alpha_2)\in\mathbb{R}_{++}^2$, the corresponding weighted-sum problem is given by
\begin{align*}
\min\limits_{x\in\mathbb{R}_+^3}
\quad
&
\alpha_1
\max_{u\in U}
\left\{
f_1(x,u)
\right\}
+
\alpha_2
\max_{u\in U}
\left\{
f_2(x,u)
\right\}
\tag{$\rm ER2_{\alpha}^{S}$}\label{ER2alphaS}
\\
\text{\rm s.t.} \quad \,
&
g_j(x,u)\leq0,
\
j=1,2,3,\
\forall u\in U.
\end{align*}

The corresponding moment-SDP formulation is obtained from problem~\eqref{ER2alphaStar} by setting $v=0$ and $\theta=0$, and then removing the redundant variables $y^0$ and $W$. 
For comparison, we use the same strictly positive weight vectors and computational setting as in case 1 to solve the resulting moment-SDP problems. 
The numerical computation gives $x_{S}^{(k)}
=
(2,0,2)^\top,
\
k=1,\ldots,9$. The first-stage solutions obtained in Cases~1 and~2 are reported in
Table~\ref{tab:numerical-results} for comparison.

As shown in Table~\ref{tab:numerical-results}, the recovered
first-stage solutions in Case~1 vary with the weight vector, whereas Case~2 gives the same solution for all the weights considered. This illustrates the greater flexibility provided by the affine adjustable second-stage decision rule.

Table~\ref{tab:case1-results} shows that the SDP optimal values and the corresponding recovered weighted objective values agree to four decimal places for all nine weight vectors. The gaps range from $9.7026\times10^{-9}$ to $4.9650\times10^{-8}$, while all maximum
robust constraint violations are numerically zero. These results are consistent with the exactness and solution recovery properties established in Theorem~\ref{thm:moment-sdp}.}
\end{example}

\begin{table}[t]
\caption{First-stage solutions obtained in Cases 1 and 2}
\label{tab:numerical-results}
\centering\footnotesize
\begin{tabular}{ccc}
\toprule
$\alpha$: Weight vectors
& $x$: First-stage solutions in Case 1
& $x^{\mathrm{S}}$: Solutions in Case 2 \\
\midrule
$(0.1,0.9)$ & $(2.5196,\ 0.0393,\ 1.9495)$ & $(2.0000,\ 0.0000,\ 2.0000)$ \\
$(0.2,0.8)$ & $(2.4059,\ 0.1173,\ 1.9701)$ & $(2.0000,\ 0.0000,\ 2.0000)$ \\
$(0.3,0.7)$ & $(2.3566,\ 0.1455,\ 1.9793)$ & $(2.0000,\ 0.0000,\ 2.0000)$ \\
$(0.4,0.6)$ & $(2.2831,\ 0.2146,\ 1.9757)$ & $(2.0000,\ 0.0000,\ 2.0000)$ \\
$(0.5,0.5)$ & $(2.1500,\ 0.3932,\ 1.9419)$ & $(2.0000,\ 0.0000,\ 2.0000)$ \\
$(0.6,0.4)$ & $(1.1215,\ 1.4514,\ 1.9296)$ & $(2.0000,\ 0.0000,\ 2.0000)$ \\
$(0.7,0.3)$ & $(0.7074,\ 1.8895,\ 1.9214)$ & $(2.0000,\ 0.0000,\ 2.0000)$ \\
$(0.8,0.2)$ & $(0.4423,\ 2.1834,\ 1.9097)$ & $(2.0000,\ 0.0000,\ 2.0000)$ \\
$(0.9,0.1)$ & $(0.2189,\ 2.4406,\ 1.8941)$ & $(2.0000,\ 0.0000,\ 2.0000)$ \\
\bottomrule
\end{tabular}
\end{table}

\begin{table}[t]
\caption{Numerical results for the two-stage model in Case 1}
\label{tab:case1-results}
\centering\footnotesize
\begin{tabular}{ccccc}
\toprule
$\alpha$: Weight vectors
& SDP value
& Recovered value
& Absolute gap
& Violation \\
\midrule
$(0.1,0.9)$ & $-4.1031$ & $-4.1031$ & $1.8261\times10^{-8}$ & $0$ \\
$(0.2,0.8)$ & $-3.9799$ & $-3.9799$ & $2.0823\times10^{-8}$ & $0$ \\
$(0.3,0.7)$ & $-3.8627$ & $-3.8627$ & $4.9650\times10^{-8}$ & $0$ \\
$(0.4,0.6)$ & $-3.7518$ & $-3.7518$ & $1.6520\times10^{-8}$ & $0$ \\
$(0.5,0.5)$ & $-3.6697$ & $-3.6697$ & $2.1240\times10^{-8}$ & $0$ \\
$(0.6,0.4)$ & $-3.6818$ & $-3.6818$ & $3.1334\times10^{-8}$ & $0$ \\
$(0.7,0.3)$ & $-3.7966$ & $-3.7966$ & $1.7462\times10^{-8}$ & $0$ \\
$(0.8,0.2)$ & $-3.9515$ & $-3.9515$ & $9.7026\times10^{-9}$ & $0$ \\
$(0.9,0.1)$ & $-4.1334$ & $-4.1334$ & $1.9973\times10^{-8}$ & $0$ \\
\bottomrule
\end{tabular}
\end{table}

\section{Conclusion}\label{sect:5}
This paper studied Geoffrion properly efficient solutions of affinely adjustable robust two-stage convex polynomial multiobjective optimization problems with spectrahedral uncertainty.
We established necessary and sufficient optimality conditions for convex polynomial data and derived an exact SOS characterization under SOS-convexity.
We further developed a finite-dimensional moment-SDP formulation based on strictly positive weighted-sum scalarization and established its optimal-value exactness and solution-recovery properties under the strict feasibility condition.
A numerical example illustrated the accuracy of solution recovery and the benefits of the adjustable two-stage model compared with its single-stage counterpart.

\subsection*{Acknowledgments.}
This work was supported by the Natural Science Foundation of Jilin Province, China, under Grant
YDZJ202401638ZYTS.
The authors would like to thank Thai Doan Chuong and Liguo Jiao for many helpful discussions.

\bibliographystyle{abbrvurl}
\bibliography{BDH}

\end{document}